\documentclass[12pt]{article}%
\usepackage{amsmath}
\usepackage{amsfonts}
\usepackage[twocolumn=false]{geometry}
\usepackage{graphicx}%
\usepackage{amssymb}
\providecommand{\U}[1]{\protect\rule{.1in}{.1in}}

\begin{document}

\title{Some Aspect of Certain two Subclass of Analytic Functions with Negative
Coefficients Defined by Rafid Operator}
\author{Alaa H. El-Qadeem \ and \ Sameerah K. Al-ghazal
\and Department of Mathematics, Faculty of Science, Zagazig University,\\Zagazig 44519, Egypt
\and ahhassan@science.zu.edu.eg\ \ \& \ \ \ z1z10z78z@gmail.com}
\date{}
\maketitle

\begin{abstract}
In this paper, we recall the subclasses $R_{\mu,p}^{\delta}(\alpha;A,B)$ and
$P_{\mu,p}^{\delta}(\alpha;A,B)$of analytic functions in the open unit disc.
Then the neighborhood properties, integral means inequalities and some results
concerning the partial sums of the functions were discussed$.$

\end{abstract}



\section*{1. Introduction}

Let $T(p)$ be the class of all p-valent functions of the from
\begin{equation}
f(z)=z^{p}-%
{\textstyle\sum\limits_{k=p+1}^{\infty}}
a_{k}z^{k}\text{ \ }(a_{k}\geq0,p\in%
\mathbb{N}
=\{1,2,3...\}),\tag{1.1}%
\end{equation}
which are analytic in the open unit disc $U=\{z\in%
\mathbb{C}
:|z|$ $<1\}$. A function $f\in T(p)$ is called p-valent starlike of order
$\alpha(0\leq\alpha<p),$ if and only if%

\begin{equation}
\operatorname*{Re}\left\{  \frac{zf^{^{\prime}}(z)}{f(z)}\right\}
>\alpha\text{ \ \ }\left(  0\leq\alpha<p;z\in U\right)  , \tag{1.2}%
\end{equation}
\noindent we denote by $T^{\ast}(p,\alpha)$ the class of all p-valent starlike
functions of order $\alpha$. Also a function $f\in T(p)$ is called p-valent
convex of order $\alpha(0\leq\alpha<p),$ if and only if%
\begin{equation}
\operatorname*{Re}\left\{  1+\frac{zf^{^{\prime\prime}}(z)}{f^{^{\prime}}%
(z)}\right\}  >\alpha\text{ \ \ }\left(  0\leq\alpha<p;z\in U\right)  ,
\tag{1.3}%
\end{equation}
we denote by $C(p,\alpha)$ the class of all p-valenty convex functions of
order $\alpha$. For more informations about the subclasses $T^{\ast}%
(p,\alpha)$ and $C(p,\alpha)$, see \cite{Silverman2}.\newline Motivated by
Atshan and Rafid see \cite{sp 1}, we introduce the following p-valent analogue
$R_{\mu,p}^{\delta}:T(p)\longrightarrow T(p):$\newline For $0\leq\mu<1$ and
$0\leq\delta\leq1,$%
\begin{equation}
R_{\mu,p}^{\delta}f(z)=\frac{1}{\Gamma(p+\delta)(1-\mu)^{p+\delta}}%
{\displaystyle\int\limits_{0}^{\infty}}
t^{\delta-1}e^{-\left(  \frac{t}{1-\mu}\right)  }f(zt)dt \tag{1.4}%
\end{equation}
Then it is easily to deduce the series representation of the function
$R_{\mu,p}^{\delta}f(z)$ as following:%

\begin{align*}
R_{\mu,p}^{\delta}(z)  &  =z^{p}-%
{\textstyle\sum\limits_{k=p+1}^{\infty}}
\frac{\Gamma(k+\delta)(1-\mu)^{k-p}}{\Gamma(p+\delta)}a_{k}z^{k},\\[1pt]
(0  &  \leq\mu<1;0\leq\delta\leq1)
\end{align*}

\noindent where $\Gamma$\ stands for Euler's Gamma function (which is valid
for all complex numbers except the non-positive integers). More operators on
the spaces of functions, see \cite{Ashwah and hassan}, \cite{4} and
\cite{Raina and Sharma}.

\noindent Let $f\ $and\ $g$\ be analytic in $U$. Then we say that the function
$g$\ is subordinate to$\ f$\ if there exists an analytic function\ in $U$ such
that\ $\left\vert w(z)\right\vert <1$ $(\forall z\in U)$\ and\ $g(z)=f(w(z))$.
For this subordination, the symbol\ $g(z)\prec f(z)$ is used. In case $f(z)$
is univalent in $U$, the subordination $g(z)\prec f(z)$ is equivalent
to\ $g(0)=f(0)$ and $\ g(U)\subset f(U)$ (see Miller and Mocanu \cite{Miller
and Mocano}).\newline For $-1\leq B<A\leq1$ and $0\leq\alpha<p$, let
$R_{\mu,p}^{\delta}(\alpha;A,B)$\ be the subclass of functions $f\in T(p)$ for which:%

\begin{equation}
\frac{z(R_{\mu,p}^{\delta}f(z))^{^{\prime}}}{R_{\mu,p}^{\delta}f(z)}%
\prec(p-\alpha)\frac{1+Az}{1+Bz}+\alpha, \tag{1.5}%
\end{equation}
that is, that
\[
R_{\mu,p}^{\delta}(\alpha;A,B)=\left\{  f\in T(p):\left\vert \tfrac
{\tfrac{z(R_{\mu,p}^{\delta}f(z))^{\prime}}{R_{\mu,p}^{\delta}f(z)}-p}%
{B\tfrac{z(R_{\mu,P}^{\delta}f(z))^{\prime}}{R_{\mu,p}^{\delta}f(z)}-\left[
Bp+(A-B)(p-\alpha)\right]  }\right\vert <1,z\in U\right\}  .
\]
Note that $\operatorname{Re}\left\{  (p-\alpha)\frac{1+Az}{1+Bz}%
+\alpha\right\}  >\tfrac{1-A+\alpha(A-B)}{1-B}.$\newline Also, for $-1\leq
B<A\leq1$ and $0\leq\alpha<p$, let $P_{\mu,p}^{\delta}(\alpha;A,B)$ be the
subclass of functions $f$ $\in T$ $(p)$ for which:%
\begin{equation}
1+\frac{z(R_{\mu,p}^{\delta}f(z))^{^{\prime\prime}}}{(R_{\mu,p}^{\delta
}f(z))^{^{\prime}}}\prec(p-\alpha)\frac{1+Az}{1+Bz}+\alpha, \tag{1.6}%
\end{equation}
For (1.5) and (1.6) it is clear that%
\begin{equation}
f(z)\in P_{\mu,p}^{\delta}(\alpha;A,B)\iff\frac{zf^{^{\prime}}(z)}{p}\in
R_{\mu,p}^{\delta}(\alpha;A,B) \tag{1.7}%
\end{equation}
The object of the present paper is to investigate the coefficients bounds,
neighborhood properties, integral means inequalities and some results
concerning partial sums for functions belonging to the subclasses $R_{\mu
,p}^{\delta}(\alpha;A,B)$ and $P_{\mu,p}^{\delta}(\alpha;A,B).$

\section*{2. Neighborhood Results}

We assume in the reminder of this paper that, $0\leq\alpha<p,$ $0\leq$
$\mu<1,0\leq\delta\leq1,-1\leq B<A\leq1,$ $p\in%
\mathbb{N}
$ and $z\in U.$ Also, we shall need the following two lemmas.\newline%
\textbf{Lemma 1 }(see \cite{Alaa Hassan}). Let the function $f(z)$ be given by
(1.1). Then $f\in R_{\mu,p}^{\delta}(\alpha;A,B),$ if and only if%

\begin{equation}%
{\textstyle\sum\limits_{k=p+1}^{\infty}}
\left[  (1-B)(k-p)+(A-B)(p-\alpha)\right]  (1-\mu)^{k-p}\frac{\Gamma
(k+\delta)}{\Gamma(p+\delta)}a_{k}\leq(A-B)(p-\alpha). \tag{2.1}%
\end{equation}
\textbf{Lemma 2 }(see \cite{Alaa Hassan}). Let the function $f(z)$ be given by
(1.1). Then $f\in P_{\mu,p}^{\delta}(\alpha;A,B),$ if and only if%

\begin{equation}%
{\textstyle\sum\limits_{k=p+1}^{\infty}}
k[(1-B)(k-p)+(A-B)(p-\alpha)](1-\mu)^{k-p}\frac{\Gamma(k+\delta)}%
{\Gamma(p+\delta)}a_{k}\leq p(A-B)(p-\alpha).\tag{2.2}%
\end{equation}
\newline Following the earlier investigations of Goodman \cite{Goodman} and
Ruscheweyh \cite{Ruscheweyh}, we recall the $\epsilon-$ neighborhood of a
function $f$ of the form (1.1) as following:\noindent%
\begin{equation}
N_{\epsilon}(f)=\left\{  g\in T(p):g(z)=z^{p}-%
{\textstyle\sum\limits_{k=p+1}^{\infty}}
b_{k}z^{k},%
{\textstyle\sum\limits_{k=p+1}^{\infty}}
k\left\vert a_{k}-b_{k}\right\vert \leq\epsilon\right\}  .\tag{2.3}%
\end{equation}
For the identity function $e(z)=z$, we immediately have%
\begin{equation}
N_{\epsilon}(e)=\left\{  g\in T(P):g(z)=z^{p}-%
{\textstyle\sum\limits_{k=p+1}^{\infty}}
b_{k}z^{k},%
{\textstyle\sum\limits_{k=p+1}^{\infty}}
kb_{k}\leq\epsilon\right\}  .\tag{2.4}%
\end{equation}
\textbf{Theorem 1}. \textit{If the function }$f(z)$\textit{\ defined by (1.1)
is in the class }$R_{\mu,p}^{\delta}(\alpha;A,B)$, \textit{then }$R_{\mu
,p}^{\delta}(\alpha;A,B)\subseteq N_{\epsilon}(e),$ where%
\[
\epsilon=\frac{(p+1)(A-B)(p-\alpha)}{\left[  (1-B)+(A-B)(p-\alpha)\right]
(1-\mu)(p+\delta)}.
\]
\textbf{Proof}. Since $f\in R_{\mu,p}^{\delta}(\alpha;A,B),$ by using Lemma 1,
we find%
\[
\frac{\left[  (1-B)+(A-B)(p-\alpha)\right]  (1-\mu)\frac{\Gamma(p+\delta
+1)}{\Gamma(p+\delta)}}{(p+1)}%
{\textstyle\sum\limits_{k=p+1}^{\infty}}
ka_{k}\leq
\]

\[%
{\textstyle\sum\limits_{k=p+1}^{\infty}}
\left[  (1-B)(k-p)+(A-B)(p-\alpha)\right]  (1-\mu)^{k-p}\frac{\Gamma
(k+\delta)}{\Gamma(p+\delta)}a_{k}\leq(A-B)(p-\alpha).
\]
Then, it is clear that%
\[%
{\textstyle\sum\limits_{k=p+1}^{\infty}}
ka_{k}\leq\frac{(p+1)(A-B)(p-\alpha)}{\left[  (1-B)+(A-B)(p-\alpha)\right]
(1-\mu)(p+\delta)}=\epsilon.
\]
This completes the proof.\newline\bigskip\textbf{Theorem 2.} \textit{If}%
\[
\epsilon=\frac{\Gamma(p+\delta)p(p+1)(A-B)(p-\alpha)}{\left[
(1-B)+(A-B)(p-\alpha)\right]  (1-\mu)\Gamma(p+\delta+1)},
\]
then $P_{\mu,p}^{\text{ }\delta}(\alpha;A,B)\subseteq N_{\epsilon}(e).$

\noindent\textbf{Proof.} For function $f\in P_{\mu,p}^{\text{ }\delta}%
(\alpha;A,B),$ of the form (1.1), from Lemma 2, we find%
\begin{align*}
\lbrack(1-B)+(A-B)(p-\alpha)](1-\mu)\frac{\Gamma(p+\delta+1)}{\Gamma
(p+\delta)}%
{\textstyle\sum\limits_{k=p+1}^{\infty}}
ka_{k} &  \leq\\%
{\textstyle\sum\limits_{k=p+1}^{\infty}}
[(k-p)(1-B)+(A-B)(p-\alpha)](1-\mu)^{k-p}\frac{\Gamma(k+\delta)}%
{\Gamma(p+\delta)}ka_{k} &  \leq p(A-B)(p-\alpha),
\end{align*}
then%
\[%
{\textstyle\sum\limits_{k=p+1}^{\infty}}
ka_{k}\leq\frac{p(A-B)(p-\alpha)}{\left[  (1-B)+(A-B)(p-\alpha)\right]
(1-\mu)(p+\delta)}=\epsilon,
\]
and the proof is completed.\newline Moreover, we will determine the
neighborhood properties for each of the following (slightly modified) function
classes $R_{\mu,p}^{\delta,\rho}(\alpha;A,B)$ and $P_{\mu,p}^{\delta,\rho
}(\alpha;A,B).$\newline A functions $f\in T(p)$ is said to be in the class
$R_{\mu,p}^{\delta,\rho}(\alpha;A,B)$\ if there exists a function $g\in
R_{\mu,p}^{\delta}(\alpha;A,B)$ such that%

\begin{equation}
\left\vert \frac{f(z)}{g(z)}-1\right\vert <1-\rho\text{ \ \ \ }(z\in
U;0\leq\rho<1).\tag{2.5}%
\end{equation}
Analogously, a function $f\in T(p)$ is said to be in the class $P_{\mu
,p}^{\delta,\rho}(\alpha;A,B)$\ if there exists a function $g\in P_{\mu
,p}^{\delta}(\alpha;A,B)$ such that the inequality (2.5) holds true.\newline
Now, using the same technique of Altintas et al. \cite{ALTINTAS}, the
neighborhood properties of the subclasses $R_{\mu,p}^{\delta,\rho}%
(\alpha;A,B)$ and $P_{\mu,p}^{\delta,\rho}(\alpha;A,B)$ are given.

\noindent\textbf{Theorem 3}.\textit{\ Let }$g\in R_{\mu,p}^{\delta}%
(\alpha;A,B)$\textit{, Suppose also that }%

\begin{equation}
\rho_{1}=1-\tfrac{\epsilon\lbrack(1-B)+(A-B)(p-\alpha)](1-\mu)(p+\delta
)}{(p+1)\{[(1-B)+(A-B)(p-\alpha)](1-\mu)(p+\delta)-(A-B)(p-\alpha
)\}},\tag{2.6}%
\end{equation}
\textit{then }%

\[
N_{\epsilon}(g)\subset R_{\mu,p}^{\delta,\rho_{1}}(\alpha;A,B).
\]
\textbf{Proof.} Assume that $f\in N_{\epsilon}(g)$. Then we find from (2.3) we get%

\[%
{\textstyle\sum\limits_{k=p+1}^{\infty}}
k\left\vert a_{k}-b_{k}\right\vert \leq\epsilon,
\]
since $g\in R_{\mu,p}^{\delta}(\alpha;A,B),$ then we have%

\[%
{\textstyle\sum\limits_{k=p+1}^{\infty}}
b_{k}\leq\frac{(A-B)(p-\alpha)}{[(1-B)+(A-B)(p-\alpha)](1-\mu)(p+\delta)},
\]
so that
\[
\left\vert \frac{f(z)}{g(z)}-1\right\vert \leq\frac{%
{\textstyle\sum\limits_{k=p+1}^{\infty}}
\left\vert a_{k}-b_{k}\right\vert }{1-%
{\textstyle\sum\limits_{k=p+1}^{\infty}}
b_{k}}\leq\frac{\epsilon}{(p+1)}\tfrac{[(1-B)+(A-B)(p-\alpha)](1-\mu
)(p+\delta)}{[(1-B)+(A-B)(p-\alpha)](1-\mu)(p+\delta)-(A-B)(p-\alpha)}%
=1-\rho_{1},
\]
provided that $\rho_{1}$ is given precisely by (2.6). Thus, by
definition,$\ f\in R_{\mu,p}^{\delta,\rho_{1}}(\alpha;A,B)$ for $\rho_{1}$
given by (2.6). This evidently completes the proof of Theorem 3.\newline
Another result regarding the subclass $P_{\mu,p}^{\delta,\rho}(\alpha;A,B)$ is
given below and the proof is omitted.

\noindent\textbf{Theorem 4}.\textit{\ If }$g\in P_{\mu,p}^{\delta}%
(\alpha;A,B)$\textit{\ and}%

\begin{equation}
\rho_{2}=1-\tfrac{\epsilon\lbrack(1-B)+(A-B)(p-\alpha)](1-\mu)(p+\delta
)}{[(1-B)+(A-B)(p-\alpha)](1-\mu)(p+\delta)(p+1)-(A-B)(p-\alpha)},\tag{2.7}%
\end{equation}
\noindent then%

\[
N_{\epsilon}(g)\subset P_{\mu,p}^{\delta,\rho_{2}}(\alpha;A,B).
\]
Now, a third neighborhood result is discussed, for this purpose we define the
subclass $H_{\mu,p}^{\text{ }\delta}(\alpha,\varphi;A,B)$\bigskip which is
related to the main subclass $R_{\mu,p}^{\delta}(\alpha;A,B),$ as
following:\newline A function $f\in T(p)$ is said to be in the class
$H_{\mu,p}^{\text{ }\delta}(\alpha,\varphi;A,B)$ if it satisfies the following
nonhomogeneous Cauchy-Euler differential equation:%

\begin{align}
z^{\text{ }2}\frac{d^{\text{ }2}f}{dz^{\text{ }2}}+2(\varphi+1)z\frac
{df(z)}{dz}+\varphi(\varphi+1)f(z)  &  =(p+\varphi)(p+\varphi+1)g(z)\tag{2.8}%
\\
(g  &  \in R_{\mu,p}^{\delta}(\alpha;A,B);\varphi>-p;).\nonumber
\end{align}
\textbf{Theorem 5}.\textit{\ If }$f\in T(p)$\textit{\ is in the subclass
}$H_{\mu,p}^{\text{ }\delta}(\alpha,\varphi;A,B),$\textit{\ then}%

\[
H_{\mu,p}^{\text{ }\delta}(\alpha,\varphi;A,B)\subset N_{\epsilon}(g),
\]
\textit{where}%
\begin{equation}
\epsilon=\frac{(p+1)(A-B)(p-\alpha)}{[(1-B)+(A-B)(p-\alpha)](1-\mu)(p+\delta
)}\left(  \frac{2(p+\varphi+1)}{(p+\varphi+2)}\right)  \tag{2.9}%
\end{equation}
\textbf{Proof}. Suppose that $f\in H_{\mu,p}^{\text{ }\delta}(\alpha
,\varphi;A,B)$ and $f$ is given by (1.1), then from (2.8) we deduce that%

\begin{equation}
a_{k}=\frac{(p+\varphi)(p+\varphi+1)}{(k+\varphi)(k+\varphi+1)}b_{k}\text{
\ \ \ \ \ }(k\geq p+1). \tag{2.10}%
\end{equation}
Moreover,%
\begin{equation}%
{\textstyle\sum\limits_{k=p+1}^{\infty}}
k\left\vert b_{k}-a_{k}\right\vert \leq%
{\textstyle\sum\limits_{k=p+1}^{\infty}}
kb_{k}+%
{\textstyle\sum\limits_{k=p+1}^{\infty}}
ka_{k}\text{ \ \ \ \ }(a_{k}\geq0,b_{k}\geq0). \tag{2.11}%
\end{equation}
by using (2.10), then (2.11) can be rewritten as following%
\begin{equation}%
{\textstyle\sum\limits_{k=p+1}^{\infty}}
k\left\vert b_{k}-a_{k}\right\vert \leq%
{\textstyle\sum\limits_{k=p+1}^{\infty}}
kb_{k}+%
{\textstyle\sum\limits_{k=p+1}^{\infty}}
\frac{(p+\varphi)(p+\varphi+1)}{(k+\varphi)(k+\varphi+1)}kb_{k}. \tag{2.12}%
\end{equation}
Next, since $g\in R_{\mu,p}^{\delta}(\alpha;A,B),$ then assertion (2.1) of the
Lemma 1 yields%

\begin{equation}%
{\textstyle\sum\limits_{k=p+1}^{\infty}}
kb_{k}\leq\frac{(p+1)(A-B)(p-\alpha)}{\left[  (1-B)+(A-B)(p-\alpha)\right]
(1-\mu)(p+\delta)}. \tag{2.13}%
\end{equation}
Finally, by making use of (2.13) on the right-hand side of (2.12), we find that%

\begin{align*}%
{\textstyle\sum\limits_{k=2}^{\infty}}
k\left\vert b_{k}-a_{k}\right\vert  &  \leq\frac{(p+1)(A-B)(p-\alpha)}{\left[
(1-B)+(A-B)(p-\alpha)\right]  (1-\mu)(p+\delta)}\left(  1+\frac{(p+\varphi
)}{(p+\varphi+2)}\right)  \\
&  =\frac{(p+1)(A-B)(p-\alpha)}{\left[  (1-B)+(A-B)(p-\alpha)\right]
(1-\mu)(p+\delta)}\left(  \frac{2(p+\varphi+1)}{(p+\varphi+2)}\right)
=\epsilon
\end{align*}
Thus, $f\in N_{\epsilon}(g).$ This, evidently, completes the proof of Theorem
5.\newline A similar result regarding the class $P_{\mu,p}^{\text{ }\delta
}(\alpha;A,B)$ can be achieved using the same techniques as performed in
Theorem 5, thus it is omitted.

\subsection*{3. Integral Means Inequalities}

In this section, we shall need the subordination lemma of Littlewood
\cite{LITTLEWOOD}.

\noindent\textbf{Lemma 3 (} \cite{LITTLEWOOD}). \textit{If the functions }$f
$\textit{\ and }$g$\textit{\ are analytic in }$U$\textit{\ with }$g(z)\prec
f(z)$\textit{\ then}%

\begin{equation}
\int_{0}^{2\pi}\left\vert g(re^{i\text{ }\theta})\right\vert ^{\tau}%
d\theta\leq\int_{0}^{2\pi}\left\vert f(re^{i\text{ }\theta})\right\vert
^{\tau}d\theta\text{ \ \ \ }(\tau>0;0<r<1). \tag{3.1}%
\end{equation}
\textbf{Theorem 6.} \textit{Let }$f\in R_{\mu,p}^{\delta}(\alpha
;A,B),$\textit{\ and suppose that }%

\begin{equation}
f_{p+1}(z)=z^{p}-\frac{(p+1)(A-B)(p-\alpha)}{[(1-B)+(A-B)(p-\alpha
)](1-\mu)(p+\delta)}z^{p+1}, \tag{3.2}%
\end{equation}
then for we have%
\begin{align}
\int_{0}^{2\pi}\left\vert f(re^{i\text{ }\theta})\right\vert ^{\tau}d\theta &
\leq\int_{0}^{2\pi}\left\vert f_{p+1}(re^{i\text{ }\theta})\right\vert ^{\tau
}d\theta\text{ }\tag{3.3}\\
(\tau &  >0,z=re^{i\text{ }\theta}(0<r<1))\nonumber
\end{align}
\textbf{Proof}. From lemma 3, it would suffice to show that%

\[
1-%
{\textstyle\sum\limits_{k=p+1}^{\infty}}
a_{k}z^{k-p}\prec1-\frac{(A-B)(p-\alpha)}{[(1-B)+(A-B)(p-\alpha)](1-\mu
)(p+\delta)}z.
\]
By setting%

\[
1-%
{\textstyle\sum\limits_{k=p+1}^{\infty}}
a_{k}z^{k-p}=1-\frac{(A-B)(p-\alpha)}{[(1-B)+(A-B)(p-\alpha)](1-\mu
)(p+\delta)}w(z).
\]
Then we find that
\begin{align*}
\left\vert w(z)\right\vert  &  =\left\vert
{\textstyle\sum\limits_{k=p+1}^{\infty}}
\frac{[(1-B)+(A-B)(p-\alpha)](1-\mu)(p+\delta)}{(A-B)(p-\alpha)}a_{k}%
z^{k-p}\right\vert \\
&  \leq\left\vert z\right\vert
{\textstyle\sum\limits_{k=p+1}^{\infty}}
\frac{[(1-B)+(A-B)(p-\alpha)](1-\mu)(p+\delta)}{(A-B)(p-\alpha)}a_{k}\\
&  \leq\left\vert z\right\vert \leq1,
\end{align*}
by using (2.1). Hence $f(z)\prec f_{p+1}(z)$ which readily yields the integral
means inequality (3.3).

\subsection*{4. Partial Sums}

In this section we will study the ratio of a function of the form (1.1) to its
sequence of partial sums defined by$f_{1}(z)=z$ and $f_{m}(z)=z^{p}-%
{\textstyle\sum\nolimits_{k=p+1}^{m}}
a_{k}z^{k}$ when the coefficients of $f(z)$ are satisfy the condition (2.1).
We will determine sharp lower bounds of $\operatorname{Re}\left(  \frac
{f(z)}{f_{m}(z)}\right)  ,$ $\operatorname{Re}\left(  \frac{f_{m}(z)}%
{f(z)}\right)  ,$ $\operatorname{Re}\left(  \frac{f^{^{\prime}}(z)}%
{f_{m}^{^{\prime}}(z)}\right)  $ and $\operatorname{Re}\left(  \frac
{f_{m}^{^{\prime}}(z)}{f^{^{\prime}}(z)}\right)  .$

\noindent In what follows, we will use the well-known result%

\[
\operatorname{Re}\left(  \frac{1-w(z)}{1+w(z)}\right)  >0\text{ \ \ \ \ }%
\left(  z\in U\right)  ,
\]
if and only if%
\[
w(z)=%
{\textstyle\sum\limits_{k=1}^{\infty}}
{}c_{k}z^{k},
\]
\newline satisfies the inequality $\left\vert w(z)\right\vert \leq\left\vert
z\right\vert .$\newline\textbf{Theorem 7.} \textit{Let }$f\in R_{\mu
,p}^{\delta}(\alpha;A,B)$\textit{, then}%

\begin{equation}
\operatorname{Re}\left(  \frac{f(z)}{f_{m}(z)}\right)  \geq1-\frac{1}{c_{m+1}%
}\text{ \ \ \ \ }\left(  z\in U,m\in%
\mathbb{N}
\right)  , \tag{4.1}%
\end{equation}
\textit{and}%

\begin{equation}
\operatorname{Re}\left(  \frac{f_{m}(z)}{f(z)}\right)  \geq\frac{c_{m+1}%
}{1+c_{m+1}}\text{ \ \ \ \ }\left(  z\in U,m\in%
\mathbb{N}
\right)  , \tag{4.2}%
\end{equation}
\textit{where}%
\begin{equation}
c_{k}=\frac{[(1-B)(k-p)+(A-B)(p-\alpha)](1-\mu)\Gamma(p+\delta+1)}%
{\Gamma(p+\delta)(A-B)(p-\alpha)}. \tag{4.3}%
\end{equation}
\textit{The estimates in (4.1) and (4.2) are sharp.}

\noindent\textbf{Proof}. Employing the same technique used by Silverman
\cite{Silverman}. The function $f\in R_{\mu,p}^{\delta}(\alpha;A,B)$ if and
only if $\sum_{k=1}^{\infty}c_{k}z^{k}\leq1.$ It is easy to verify that
$c_{k+1}>c_{k}>1.$ Thus%

\begin{equation}%
{\textstyle\sum\limits_{k=p}^{m}}
a_{k}+c_{m+1}%
{\textstyle\sum\limits_{k=m+p}^{\infty}}
a_{k}\leq%
{\textstyle\sum\limits_{k=p+1}^{\infty}}
c_{k}a_{k}<1. \tag{4.4}%
\end{equation}
Now, setting
\[
c_{m+1}\left\{  \frac{f(z)}{f_{m}(z)}-\left(  1-\frac{1}{c_{m+1}}\right)
\right\}  =\frac{1-%
{\textstyle\sum\limits_{k=p}^{m}}
a_{k}z^{k-p}-c_{m+1}%
{\textstyle\sum\limits_{k=m+p}^{\infty}}
a_{k}z^{k-p}}{1-%
{\textstyle\sum\limits_{k=p}^{m}}
a_{k}z^{k-p}}=\frac{1+D(z)}{1+E(z)},
\]
and $\frac{1+D(z)}{1+E(z)}=\frac{1-w(z)}{1+w(z)},$ then we have%

\[
w(z)=\frac{E(z)-D(z)}{2+D(z)+E(z)}=\frac{c_{m+1}%
{\textstyle\sum\limits_{k=m+p}^{\infty}}
a_{k}z^{k-p}}{2-2%
{\textstyle\sum\limits_{k=p}^{m}}
a_{k}z^{k-p}-c_{m+1}%
{\textstyle\sum\limits_{k=m+p}^{\infty}}
a_{k}z^{k-p}}
\]
which implies
\[
\left\vert w(z)\right\vert \leq\frac{c_{m+1}%
{\textstyle\sum\limits_{k=m+p}^{\infty}}
a_{k}}{2-2%
{\textstyle\sum\limits_{k=p}^{m}}
a_{k}-c_{m+1}%
{\textstyle\sum\limits_{k=m+p}^{\infty}}
a_{k}}.
\]
Hence $\left\vert w(z)\right\vert \leq1,$ if and only if%
\[%
{\textstyle\sum\limits_{k=p}^{m}}
a_{k}+c_{m+1}%
{\textstyle\sum\limits_{k=m+p}^{\infty}}
a_{k}\leq1
\]
which is true by (4.4). This readily yields (4.1).\newline Now consider the function%

\begin{equation}
f(z)=1-\frac{z^{m+1}}{c_{m+1}} \tag{4.5}%
\end{equation}
Thus $\frac{f(z)}{f_{m}(z)}=1-\frac{z^{m}}{c_{m+1}}.$ Letting
$z\longrightarrow1^{-},$ then $f(z)=1-\frac{1}{c_{m+1}}.$ So $f(z)$ given by
(4.5) satisfies the sharp result in (4.1). shows that the bounds in (4.1) are
best possible for each $m\in%
\mathbb{N}
.$

\noindent Similarly, setting%
\[
\left(  1+c_{m+1}\right)  \left\{  \frac{f_{m}(z)}{f(z)}-\frac{c_{m+1}%
}{1+c_{m+1}}\right\}  =\frac{1-%
{\textstyle\sum\limits_{k=p}^{m}}
a_{k}z^{k-p}+c_{m+1}%
{\textstyle\sum\limits_{k=m+p}^{\infty}}
a_{k}z^{k-p}}{1-%
{\textstyle\sum\limits_{k=p}^{m}}
a_{k}z^{k-p}}=\frac{1-w(z)}{1+w(z)},
\]
where%
\[
\left\vert w(z)\right\vert \leq\frac{(1+c_{m+1})%
{\textstyle\sum\limits_{k=m+p}^{\infty}}
a_{k}}{2-2%
{\textstyle\sum\limits_{k=p}^{m}}
a_{k}+(1-c_{m+1})%
{\textstyle\sum\limits_{k=m+p}^{\infty}}
a_{k}}.
\]
Now $\left\vert w(z)\right\vert \leq1,$ if and only if%

\[%
{\textstyle\sum\limits_{k=p}^{m}}
a_{k}+c_{m+1}%
{\textstyle\sum\limits_{k=m+p}^{\infty}}
a_{k}\leq1,
\]
which readily implies the assertion (4.2). The estimate in (4.2) is sharp with
the extremal function $f(z)$ given by (4.5). This completes the proof of the
theorem 7.

\noindent Following similar steps to that followed in Theorem 7, we can state
the following theorem

\noindent\textbf{Theorem 8.} \textit{Let }$f\in R_{\mu,p}^{\delta}%
(\alpha;A,B)$\textit{, then}%

\begin{equation}
\operatorname{Re}\left(  \frac{f^{\text{ }^{\prime}}(z)}{f_{m}^{^{\text{
}\prime}}(z)}\right)  \geq1-\frac{m+1}{c_{m+1}}\text{ \ \ \ \ }\left(  z\in
U,m\in%
\mathbb{N}
\right)  , \tag{4.6}%
\end{equation}
\textit{and}%

\begin{equation}
\operatorname{Re}\left(  \frac{f_{m}^{^{\text{ }\prime}}(z)}{f^{\text{
}^{\prime}}(z)}\right)  \geq\frac{c_{m+1}}{m+1+c_{m+1}}\text{ \ \ \ \ }\left(
z\in U,m\in%
\mathbb{N}
\right)  , \tag{4.7}%
\end{equation}
In both cases, the extremal function$f(z)$ is as defined in (4.5).\newline%
\textbf{Proof}. We prove only (4.6), which is similar in spirit to the proof
of theorem 7 and similariy we proof (4.7). We write

\textbf{\ }%
\[
\frac{c_{m+1}}{m+1}\left\{  \frac{f^{\ ^{\prime}}(z)}{f_{m}^{\text{ \ }\prime
}(z)}-\left(  1-\frac{m+1}{c_{m+1}}\right)  \right\}  =\frac{1-%
{\textstyle\sum\limits_{k=p}^{m}}
a_{k}z^{k-1}-\frac{c_{m+1}}{m+1}%
{\textstyle\sum\limits_{k=m+p}^{\infty}}
a_{k}z^{k-1}}{1-%
{\textstyle\sum\limits_{k=p}^{m}}
a_{k}z^{k-1}}=\frac{1+D(z)}{1+E(z)}
\]
and $\frac{1+D(z)}{1+E(z)}=\frac{1-w(z)}{1+w(z)},$ then we have%

\[
w(z)=\frac{E(z)-D(z)}{2+D(z)+E(z)}=\frac{\frac{c_{m+1}}{m+1}%
{\textstyle\sum\limits_{k=m+p}^{\infty}}
a_{k}z^{k-1}}{2-2%
{\textstyle\sum\limits_{k=p}^{m}}
a_{k}z^{k-1}-\frac{c_{m+1}}{m+1}%
{\textstyle\sum\limits_{k=m+p}^{\infty}}
a_{k}z^{k-1}}
\]
which implies
\[
\left\vert w(z)\right\vert \leq\frac{\frac{c_{m+1}}{m+1}%
{\textstyle\sum\limits_{k=m+p}^{\infty}}
a_{k}}{2-2%
{\textstyle\sum\limits_{k=p}^{m}}
a_{k}-\frac{c_{m+1}}{m+1}%
{\textstyle\sum\limits_{k=m+p}^{\infty}}
a_{k}}.
\]
Hence $\left\vert w(z)\right\vert \leq1,$ if and only if%
\[%
{\textstyle\sum\limits_{k=p}^{m}}
a_{k}+\frac{c_{m+1}}{m+1}%
{\textstyle\sum\limits_{k=m+p}^{\infty}}
a_{k}\leq1.
\]

\end{document}